\documentclass[11pt]{article}
\usepackage{latexsym}
\usepackage{amsfonts}
\usepackage{amssymb}
\usepackage{amscd}
\usepackage{array}
\usepackage{amsmath}
\usepackage{epsfig}
\usepackage{longtable}
\usepackage{a4}
\begin{document}
\title{\textbf{Stringy Hodge numbers of strictly canonical nondegenerate singularities}\\ \date{}}
\author{Jan Schepers\footnote{Supported by VICI grant 639.033.402 from the Netherlands Organisation for
Scientific Research (NWO). During the completion of this paper, the author was a Postdoctoral Fellow of the Research Foundation - Flanders (FWO).\newline
\emph{Address}: Jan Schepers, Katholieke Universiteit Leuven, Department of Mathematics, Celes\-tij\-nenlaan 200B, 3001 Leuven, Belgium.\newline
\emph{E-mail}: \texttt{janschepers1@gmail.com}.}}
\maketitle

\vspace{-0.7cm}

\begin{center}
\textit{Dedicated to Joost van Hamel}
\end{center}

\vspace{0.3cm}

\begin{abstract}
\noindent We describe a class of isolated nondegenerate
hypersurface singularities that give a polynomial contribution to
Batyrev's stringy $E$-function. These singularities are obtained by
imposing a natural condition on the facets of the Newton
polyhedron, and they are strictly canonical. We prove that
Batyrev's conjecture concerning the nonnegativity of stringy Hodge
numbers is true for complete varieties with such singularities, under some additional hypotheses on the defining polynomials (e.g.\ convenient or weighted homogeneous).
The proof uses combinatorics on lattice polytopes. The results
form a strong generalisation of previously obtained results for
Bries\-korn singularities.
\end{abstract}

\section{Introduction}
\noindent \textbf{1.1.} Batyrev defined now more than a decade ago the
\textit{stringy $E$-function} for complex algebraic varieties with log
terminal singularities \cite{Batyrev}. It is a rational function in two
variables $u,v$ if the singularities are Gorenstein and canonical. Batyrev made
moreover the following fascinating conjecture: if the stringy $E$-function of a
projective variety is a polynomial $\sum_{p,q} a_{p,q}u^pv^q$ then the \textit{stringy Hodge numbers}
$h_{st}^{p,q}:=(-1)^{p+q}a_{p,q}$ are nonnegative. This conjecture is motivated by the fact that stringy
Hodge numbers share many other properties with usual Hodge numbers of smooth projective varieties.\\

\noindent \textbf{1.2.} It is known that the stringy $E$-function of toric varieties with Gorenstein
(and hence canonical) singularities is a polynomial and that Batyrev's conjecture is true for such
complete varieties (\cite{Batyrev} and \cite{MustataPayne}). The same remarks apply to varieties with Gorenstein
quotient singularities. Yasuda relates the stringy Hodge numbers in that case with orbifold cohomology \cite{Yasuda}.
Together with Veys we proved in \cite{SchepersVeys} Batyrev's conjecture in full generality for threefolds and also
for a class of isolated singularities in dimension $\geq 4$ (see Theorem 3.9 for the precise statement). The
disadvantage of that theorem in higher dimension is that the conditions demanded for the singularities will prevent
the stringy $E$-function in many examples from being a polynomial.\\

\noindent \textbf{1.3.} In this paper we focus our attention on
nondegenerate hypersurface singularities. This is a very
computable class of singularities due to the connection to the
toric world and hence they often serve as a testing ground for
open problems. They are defined using the Newton polyhedron
$\Gamma(f)$ associated to their equation $f=0$. This Newton
polyhedron induces a decomposition $\Delta_f$ of the first orthant
of the dual space into cones, called the \textit{first Varchenko
subdivision} in \cite{Stepanov}. We study isolated nondegenerate
singularities with a natural extra condition on the 1-dimensional
cones of $\Delta_f$, i.e.\ the cones of $\Delta_f$ associated to
codimension 1 faces (called \textit{facets}) of $\Gamma(f)$. We
call the subdivision $\Delta_f$ then \textit{crepant} (see
Definition 4.8 and Remark 4.9 (1)). These singularities are
strictly canonical and that means that they are part of the
`worst' untreated case in the aforementioned theorem of
\cite{SchepersVeys}, but on the other hand they are in a sense the
best chance for obtaining a polynomial stringy $E$-function
(Remark 5.4 (4)). The following theorem is our main result. It
gives a strong generalisation of the results
from \cite{SchepersVeys2} about Brieskorn singularities. See Proposition 4.11 and Theorem 5.3.\\

\noindent \textbf{1.4. Theorem.} \textsl{Let $V$ be an algebraic variety whose singularities are analytically isomorphic to isolated nondegenerate singularities with crepant first Varchenko
subdivision. Then the stringy $E$-function of $V$ is a polynomial. If $V$ is complete and if the defining polynomial $f$ of each of the singularities is convenient (i.e.\ $f$ contains a nonzero term $a_ix_i^{b_i}$ for each variable $x_i$) or if the set of compact facets of $\Gamma(f)$ has a unique maximal element (e.g.\ $f$ weighted homogeneous) then 
the stringy Hodge numbers of $V$ are nonnegative.}\\

\noindent \textbf{1.5.} The proof of this theorem mainly uses combinatorics of lattice polytopes: we will see that the ingredients to compute the contribution of the singularity to the stringy $E$-function can be expressed in terms of lattice polytopes (Theorem 4.14). The
condition of having a crepant first Varchenko subdivision is crucial for this, together with the formula
for the Hodge-Deligne polynomial of a nondegenerate hypersurface in the torus from \cite{BatyrevBorisov} and
the formula for the contribution of the singularity itself from \cite{SchepersVeys2}.\\

\noindent \textbf{1.6.} This paper is organised as follows. In
Section 2 we gather all the combinatorial definitions that we need
and we prove a few useful lemma's. That section is self-contained
and can be read separately from the rest of the paper. In Section
3 we review Batyrev's definitions of the stringy $E$-function and
the stringy Hodge numbers. We recall the basic facts about
nondegenerate singularities in Section 4. There we also define
nondegenerate singularities with crepant first Varchenko
subdivision and we prove that the contribution of such
singularities to the stringy $E$-function is a polynomial.
Finally, in Section 5 we prove the nonnegativity statement of the stringy Hodge numbers from Theorem 1.4.\\

\noindent \textbf{Acknowledgements.} Part of this work was carried
out during a stay at the Institut des Hautes \'Etudes
Scientifiques (IH\'ES). I am very grateful that I was given the
possibility to work there. I also want to thank Ann Lemahieu and Wim Veys for
helpful discussions.

\section{Combinatorial preliminaries}
\noindent \textbf{2.1.} In this section we summarise some combinatorial aspects of Eulerian posets and lattice
polytopes that are used later.\\

\noindent \textbf{2.2.} Let $\mathcal{P}$ be a finite poset (i.e.\
\textbf{p}artially \textbf{o}rdered \textbf{set}). If $x,y \in
\mathcal{P}$ and $x\leq y$, then the interval $[x,y]$ is the set
$\{z\in \mathcal{P}\, |\, x\leq z\leq y\}$. The dual
$\mathcal{P}^*$ of $\mathcal{P}$ is obtained by taking the same
underlying set with inverted partial order relation. We assume
that $\mathcal{P}$ has a minimal element $\hat{0}$ and a maximal
element $\hat{1}$ and that every maximal chain $\hat{0}=x_0 < x_1
< \cdots < x_{d-1} <x_d = \hat{1}$ has the same length $d$. In
that case there exists a unique \textit{rank function} $\rho:
\mathcal{P} \to \{0,\ldots,d\}$ such that $\rho(x)$ equals the
length of a saturated chain in the interval $[\hat{0},x]$. One calls
$\mathcal{P}$ then \textit{graded of rank} $d$. If every
nontrivial interval in $\mathcal{P}$ has the same number of
elements of even and odd rank, then $\mathcal{P}$ is called
\textit{Eulerian}. There is an equivalent formulation in terms of
the M\"obius function. This function $\mu$ is defined on pairs
$(x,y)\in \mathcal{P}\times\mathcal{P}$ with $x\leq y$ in the
following inductive way:
\[ \mu(x,x)=1, \text{ for all } x\in \mathcal{P},\]
\[ \mu(x,y) = - \sum_{x\leq z < y} \mu(x,z), \text{ for } x<y.\]
A finite graded poset is then Eulerian if and only if $\mu(x,y)=(-1)^{\rho(y)-\rho(x)}$ for all $x\leq y$.
It is easy to see that an interval in an Eulerian poset is again Eulerian and that the dual $\mathcal{P}^*$ of an
Eulerian poset $\mathcal{P}$ is Eulerian as well.\\

\noindent \textbf{2.3. Definition.} Let $\mathcal{P}$ be an
Eulerian poset of rank $d$. Define
$g(\mathcal{P},t),h(\mathcal{P},t)\in \mathbb{Z}[t]$ by the
following recursive rules:
\[ g(\mathcal{P},t)=h(\mathcal{P},t)=1 \text{ if } d=0,\]
\[ h(\mathcal{P},t) = \sum_{\hat{0}< x \leq \hat{1}} (t-1)^{\rho(x)-1} g([x,\hat{1}],t) \text{ if } d>0, \]
\[ g(\mathcal{P},t) = \tau_{< d/2} \bigl((1-t)h(\mathcal{P},t)\bigr) \text{ if } d>0, \]
where $\tau_{<r}: \mathbb{Z}[t] \to \mathbb{Z}[t] $ is the
truncation operator defined by
\[ \tau_{<r}\biggl(\sum_i a_it^i \biggr) = \sum_{i<r} a_it^i.\]

\noindent This definition was given by Stanley \cite[\S
2]{Stanley87}. In fact the polynomials defined by Stanley for
$\mathcal{P}$, are in our notation $g(\mathcal{P}^*,t)$ and
$h(\mathcal{P}^*,t)$. See also
\cite[Def.\ 2.4]{BatyrevBorisov}. For $d>0$, $\deg h(\mathcal{P},t)=d-1$ and $\deg g(\mathcal{P},t)\leq (d-1)/2$.\\

\noindent \textbf{2.4.} These polynomials have the following
properties, with $\mathcal{P}$ an Eulerian poset of rank $d>0$:
\begin{itemize}
\item[(1)] $h(\mathcal{P},t)= t^{d-1} h(\mathcal{P},t^{-1})$,
\item[(2)] $\displaystyle{\sum_{\hat{0} \leq x\leq \hat{1}}}
g([\hat{0},x],t)\, g([x,\hat{1}]^*,t) (-1)^{\rho(\hat{1})-\rho(x)}
=0\quad $ and also

$\qquad \displaystyle{\sum_{\hat{0} \leq x\leq \hat{1}}}
(-1)^{\rho(x)-\rho(\hat{0})} g([\hat{0},x]^*,t)\, g([x,\hat{1}],t)
=0$.
\end{itemize}
The first one is proved in \cite[Thm.\ 2.4]{Stanley87} and the second one, called Stanley's convolution property,
in \cite[Cor.\ 8.3]{Stanley92}.\\


\noindent \textbf{2.5.} Let $P$ be a \textit{lattice polytope} in
$\mathbb{R}^n$ (i.e.\ the convex hull of a finite number of points
with vertices in $\mathbb{Z}^n$). The dimension of $P$ is the
dimension of the smallest affine subspace of $\mathbb{R}^n$
containing $P$. We also allow $P=\emptyset$ as a polytope of
dimension $-1$. A \textit{face} of a lattice polytope $P$ is any
intersection of $P$ with a hyperplane $H$ in $\mathbb{R}^n$ such
that $P$ is completely contained in one of the two closed
halfspaces determined by $H$. A \textit{facet} is a face of
codimension 1. The empty set and $P$ itself are also considered as
faces of $P$, but are called improper faces. For faces $F,F'$ of
$P$ we write $F\leq F'$ if $F\subseteq F'$. In this way the set of
faces of $P$ becomes an Eulerian poset with rank function
$\rho(F)=\dim F+1$ \cite[p.122]{Stanley97}. We denote this poset
by $\mathcal{P}(P)$. The dual poset $\mathcal{P}(P)^*$ is also of
the form $\mathcal{P}(Q)$ for a lattice polytope $Q$. For such
posets we have the following properties:
\begin{itemize}
\item[(1)] $g(\mathcal{P}(P),t)$ and $h(\mathcal{P}(P),t)$ have nonnegative coefficients,
\item[(2)] $g(\mathcal{P}(P),t)=1$ if and only if $P$ is a simplex.
\end{itemize}
The proof of the first statement uses the connection with toric
varieties and their intersection cohomology,
see Theorem 3.1 and Corollary 3.2 from \cite{Stanley87}.\\

\noindent \textbf{2.6.} Braden and MacPherson define relative
$g$-polynomials by proving the following
\cite[Prop.\ 2]{BradenMacPherson} for arbitrary polytopes (they do not restrict to polytopes with integer vertices).\\

\noindent \textbf{Proposition.} \textsl{There is a unique family
of polynomials $g(P,F,t)\in \mathbb{Z}[t]$ associated to a
polytope $P$ and a face $F$ of $P$, satisfying the following
relation: for all $P,F$, we have}
\[ \sum_{F\leq E\leq P} g(E,F,t)\, g([E,P]^*,t) =
g(\mathcal{P}(P)^*,t). \tag{1} \]

\noindent We remark that Braden and MacPherson use the same
definition as Stanley for $g$-polynomials, so compared to their
formula we have to use dual posets at the appropriate places. But
our notation $g(P,F,t)$ corresponds to theirs. Relative
$g$-polynomials have the following properties:
\begin{itemize}
\item[(1)] the coefficients of $g(P,F,t)$ are nonnegative if $P$ is a lattice polytope (again proved using intersection cohomology, see
\cite[Thm.\ 4]{BradenMacPherson}),
\item[(2)] $g(P,P,t)=g(\mathcal{P}(P)^*,t)$ and if
$P\neq \emptyset$ then $g(P,\emptyset,t)=0$.
\end{itemize}

\noindent \textbf{2.7.} If $G$ is a face of $P$ then it is not
hard to see that there exists a lattice polytope $P/G$ whose poset
of faces equals $[G,P]$ (by a slight adaptation of the
construction in the introduction of \cite{BradenMacPherson}). Let
$F$ be a face of $P$ that contains $G$. Then $P/G$ has a face
corresponding to $F$ that we denote by $F/G$. Formula (1) applied
for $P/G$ and $F/G$ can be written as
\[ \sum_{F\leq E\leq P} g(E/G,F/G,t)\, g([E,P]^*,t) =
g([G,P]^*,t). \] Below we will use the notation $g([G,E],[G,F],t)$
for $g(E/G,F/G,t)$.\\

\noindent \textbf{2.8.} Let $P$ be a lattice polytope of dimension
$d\geq 0$. Denote by $f_P(m)$ the number of lattice points $|mP
\cap \mathbb{Z}^n|$ for $m\in \mathbb{Z}_{>0}$. It is well known
that the so-called Ehrhart generating series $1+\sum_{m> 0}
f_P(m)\, t^m$ can be written in the form
\[  \frac{h^*_P(t)}{(1-t)^{d+1}},      \]
where $h^*_P(t)$ is a polynomial of degree $s\leq d$ with nonnegative integer coefficients.
Moreover, $l=d+1-s$ is the smallest integer such that $lP$ contains a lattice point in its relative interior. \\

\noindent \textbf{2.9.} Let $p$ be a vertex of a positive
dimensional polytope $P$. The \textit{closed star neighbourhood}
$star_{\partial P}(p)$ of $p$ in the boundary $\partial P$ of $P$
is the following set of faces:
\[ \{F \text{ face of }  P \,|\, F \text{ is a face of a proper face } Q \text{ of } P \text{ with } p\in Q\}.\]
The group $AGL(n,\mathbb{Z})$ consists of affine transformations
$A$ of $\mathbb{R}^n$ such that $A(\mathbb{Z}^n)=\mathbb{Z}^n$.
Two lattice polytopes $P$ and $Q$ in $\mathbb{R}^n$ are called
\textit{isomorphic} if there exists such an affine transformation
$A$ with $A(P)=Q$. In that case we clearly have
$h_P^*(t)=h_Q^*(t)$. Let $Q$ be a lattice polytope in
$\mathbb{R}^{n-1}$. We define the \textit{standard lattice
pyramid} $\Pi(Q)$ over $Q$ as the convex hull of $Q \times \{0\}$
and $(0,\ldots,0,1)$ in $\mathbb{R}^{n}=\mathbb{R}^{n-1}\times
\mathbb{R}$. A polytope $P$ in $\mathbb{R}^n$ is called a
\textit{lattice pyramid over a facet} $F$ of $P$ if $P$ is
isomorphic to a standard lattice pyramid $\Pi(Q)$ such that $F$
corresponds to $Q$ under this isomorphism. It is not hard to prove
that then
$h_P^*(t)=h_F^*(t)$. This generalises in the following lemma.\\

\noindent \textbf{2.10. Lemma.} \textsl{Let $P$ be a lattice
polytope in $\mathbb{R}^n$ of dimension $>0$ and $p$ a vertex of
$P$. Assume that for all facets $F$ of $P$ not containing $p$, the
convex hull of $F$ and $p$ is a lattice pyramid over $F$. Then}
\[ h_P^*(t) = \sum_{\substack{G \text{ proper face}\\ G \notin star_{\partial P}(p)}} h_G^*(t)\,
(t-1)^{\dim P -\dim G -1}.\]

\noindent \textit{Proof.} We may assume that the vertex $p$ lies
at the origin of $\mathbb{R}^n$ and hence we can consider it as a
vertex of all multiples $mP$ of $P$. First we look at the poset $\mathcal{Q} = \{ G\leq P\, |\, G\notin star_{\partial P}(p)\}$. Note that every interval in $\mathcal{Q}$ is Eulerian and hence $\mu_{\mathcal{Q}}(F,G) = (-1)^{\dim G - \dim F}$ for all $F\leq G$ in $\mathcal{Q}$. For $G\in \mathcal{Q}$ we denote the convex hull of $G$ and $p$ by $(G,p)$. We define the function
\[ g_m: \mathcal{Q} \to \mathbb{Z}_{\geq 0}: G \mapsto |(m(G,p)
\setminus \{p\}) \cap \mathbb{Z}^n| \]
and the function $f_m: \mathcal{Q} \to \mathbb{Z}_{\geq 0}$ that for $G\in \mathcal{Q}$ counts the number of integer points in $m(G,p)\setminus \{p\}$ that is not contained in any $m(F,p)$ with $F\in \mathcal{Q}$ and $F<G$. Then 
\[ g_m(G) = \sum_{\substack{F\in \mathcal{Q}\\ F\leq G}} f_m(F) \]
and hence we can apply the M\"obius inversion formula \cite[Prop.\ 3.7.1]{Stanley97} to $\mathcal{Q},f_m,g_m$ to conclude that 
\[ 0 = f_m(P) = \sum_{G \in \mathcal{Q}} (-1)^{\dim P - \dim G} |(m(G,p)
\setminus \{p\}) \cap \mathbb{Z}^n|.\] 
Thus 
\[ |(mP\setminus \{p\}) \cap \mathbb{Z}^n| = \sum_{\substack{G \text{ proper face}\\ G
\notin star_{\partial P}(p)}} (-1)^{\dim P - \dim G -1} |(m(G,p)
\setminus \{p\}) \cap \mathbb{Z}^n|. \tag{*}\]
Secondly, $\mathcal{Q}'=star_{\partial P}(p)\cup \{P\}$ is a finite graded
poset with $\emptyset = \hat{0}$ and $P=\hat{1}$, where each interval $[x,y]$ with $y\neq \hat{1}$ is
Eulerian. So by definition of the M\"obius function we have
\[ \mu_{\mathcal{Q}'}(\hat{0},\hat{1}) = - \sum_{G\in star_{\partial P}(p)} (-1)^{\dim G +1}.\]
By \cite[Prop.\ 3.8.8]{Stanley97} this equals the reduced topological Euler characteristic of the space
\[ \bigcup_{G\in star_{\partial P}(p)} G.\] This space is contractible to $\{p\}$ and hence its reduced Euler
characteristic is 0. Because $\mathcal{P}(P)$ is Eulerian, we
deduce that
\[ 1= \sum_{\substack{G \text{ proper face}\\ G \notin star_{\partial P}(p)}} (-1)^{\dim P - \dim G -1} . \tag{**}  \]
Adding $(*)$ and $(**)$ gives
\[ |mP \cap \mathbb{Z}^n| = \sum_{\substack{G \text{ proper face}\\ G \notin star_{\partial P}(p)}} (-1)^{\dim P
- \dim G -1} |m(G,p) \cap \mathbb{Z}^n|.\] So the Ehrhart series
of $P$ equals
\[ \sum_{\substack{G \text{ proper face}\\ G \notin star_{\partial P}(p)}} (-1)^{\dim P - \dim G -1} \biggl( 1 +
\sum_{m \in \mathbb{Z}_{>0}}|m(G,p) \cap \mathbb{Z}^n|\,t^m\biggr) .\]
Multiplying by $(1-t)^{\dim P +1}$ one gets
\[ h_P^*(t) = \sum_{\substack{G \text{ proper face}\\ G \notin star_{\partial P}(p)}} h_G^*(t) (-1)^{\dim P -
 \dim G -1}(1-t)^{\dim P -\dim G -1}.\]

\vspace{-0.6cm}

\hfill $\blacksquare$\\

\noindent \textbf{2.11. Definition} \cite[Def.\
5.3]{BorisovMavlyutov}\textbf{.} Let $P$ be a lattice polytope. We
define the polynomial $\widetilde{S}(P,t)\in \mathbb{Z}[t]$ by the
formula
\[ \widetilde{S}(P,t) = \sum_{\emptyset \leq F \leq P} (-1)^{\dim P - \dim F} h_F^*(t)\, g([F,P],t),\]
where we sum over all faces of $P$, with $h_F^*(t)$ the $h^*$-polynomial of the lattice polytope $F$ (with
$h_{\emptyset}^*(t)=1$) and with $g([F,P],t)$ the $g$-polynomial of the interval $[F,P]$ in the Eulerian poset
$\mathcal{P}(P)$. Note that $\widetilde{S}(\emptyset,t)=1$ and $\widetilde{S}(P,t)=0$ if $\dim P=0$.\\

\noindent \textbf{2.12.} The polynomial $\widetilde{S}(P,t)$ has
the following properties.
\begin{itemize}
\item[(1)] $\deg(\widetilde{S}(P,t)) \leq \dim P$.
\item[(2)] $\widetilde{S}(P,0)=0$ if $\dim P \geq 0$ since a $h^*$-
and a $g$-polynomial always have constant coefficient 1 and since
$P$ has an equal number of even- and odd-dimensional faces.
\item[(3)] The coefficients of $\widetilde{S}(P,t)$ are nonnegative. In \cite[Prop.\ 5.5]{BorisovMavlyutov} they are
interpreted as the dimensions of the pieces of the pure Hodge
structure on the lowest weight part of the middle cohomology of a
nondegenerate affine hypersurface in the maximal torus of the
toric variety associated to $P$.
\item[(4)] For instance from this description one has the reciprocity law
$\widetilde{S}(P,t) = t^{\dim P +1} \widetilde{S}(P,t^{-1})$
\cite[Rem.\ 5.4]{BorisovMavlyutov}.
\end{itemize}
For more information on this polynomial we refer to \cite[\S 4]{BatyrevNill}.\\

\noindent \textbf{2.13.} We will need an extension of the
definition of the $\widetilde{S}$-polynomial. An \textit{order ideal} $\mathcal{I}$ in a poset $\mathcal{P}$ is a subset for which $x\in \mathcal{I}$ and $y\leq x$ imply $y\in \mathcal{I}$.\\

\noindent \textbf{Definition.} Let $P$ be a lattice polytope and
$\mathcal{I} \varsubsetneq \mathcal{P}(P)$ be an order ideal of $\mathcal{P}(P)$. We define the polynomial
$\widetilde{S}(P,\mathcal{I},t)\in \mathbb{Z}[t]$ by the formula
\[ \widetilde{S}(P,\mathcal{I},t) = \sum_{\substack{\emptyset \leq F \leq P\\ F \notin \mathcal{I}}}
(-1)^{\dim P - \dim F} h_F^*(t)\, g([F,P],t).\]

\noindent Note that
\begin{itemize}
\item[(1)] $\deg \widetilde{S}(P,\mathcal{I},t)\leq \dim P$,
\item[(2)] $\widetilde{S}(P,\emptyset,t) = \widetilde{S}(P,t)$,
\item[(3)] $\widetilde{S}(P,\{\emptyset\},t) = \widetilde{S}(P,t)+(-1)^{\dim P}g(\mathcal{P}(P),t)$. In particular, if $P$ is odd-dimensional then the constant coefficient of $\widetilde{S}(P,\{\emptyset\},t)$ is $-1$. 
\end{itemize}
Hence in general there can be negative coefficients in $\widetilde{S}(P,\mathcal{I},t)$. For now, we study the special case where $\mathcal{I} = \{ F \leq P \,|\, F\not\geq Q\}$ for a fixed face $Q$ of $P$. We denote $\widetilde{S}(P,\mathcal{I},t)$ then by $\widetilde{S}(P,Q,t)$ and we show in Corollary 2.15 that $\widetilde{S}(P,Q,t)$ has nonnegative coefficients.\\

\noindent \textbf{2.14. Proposition.} \textsl{Let $P$ be a lattice
polytope and $Q'\leq Q$ be faces of $P$. For a face $F$ of $P$
denote by $F\vee Q$ the unique smallest face of $P$ containing $F$
and $Q$. Then
\[\widetilde{S}(P,Q,t) = \widetilde{S}(P,Q',t) +
\sum_{\substack{Q'\leq F < P\\ Q \not\leq F}} g([F,P],[F,F\vee
Q],t)\, \widetilde{S}(F,Q',t) .
\]}

\noindent \textit{Proof.} We work by induction on $\dim P - \dim
Q'$. The case $P=Q=Q'$ is trivial. So assume $\dim P - \dim Q'
>0$. By definition
\[ \widetilde{S}(P,Q,t) = \widetilde{S}(P,Q',t) +
\sum_{\substack{Q'\leq F < P\\ Q \not\leq F}} (-1)^{\dim P - \dim
F -1} h_F^*(t)\, g([F,P],t). \tag{2}
\]
By Stanley's convolution property we have
\begin{eqnarray*}
g([F,P],t) &=& \sum_{\substack{F\leq G < P\\ Q \not\leq G}}
(-1)^{\dim P - \dim G -1} g([F,G],t)\, g([G,P]^*,t)\\ & & +
\sum_{\substack{F\leq G' < P\\ Q \leq G'}} (-1)^{\dim P - \dim G'
-1} g([F,G'],t)\, g([G',P]^*,t).
\end{eqnarray*}
We put this in (2) and exchange the sums to find
\begin{eqnarray*}
\widetilde{S}(P,Q,t) &=& \widetilde{S}(P,Q',t) +
\sum_{\substack{Q'\leq G < P\\ Q
\not\leq G}} g([G,P]^*,t) \, \widetilde{S}(G,Q',t) \\
& & + \sum_{Q\leq G' < P} g([G',P]^*,t)
\bigl(\widetilde{S}(G',Q',t)- \widetilde{S}(G',Q,t) \bigr).
\end{eqnarray*}
By the induction hypothesis this becomes
\begin{eqnarray*} &=& \widetilde{S}(P,Q',t) + \sum_{\substack{Q'\leq G < P\\ Q
\not\leq G}} g([G,P]^*,t) \, \widetilde{S}(G,Q',t) \\
& & - \sum_{Q\leq G' < P} g([G',P]^*,t) \sum_{\substack{Q'\leq E <
G'\\ Q \not\leq E}} g([E,G'],[E,E\vee Q],t)\,
\widetilde{S}(E,Q',t) 
\end{eqnarray*}
\begin{eqnarray*}
&=& \widetilde{S}(P,Q',t) + \sum_{\substack{Q'\leq G < P\\ Q
\not\leq G}} g([G,P]^*,t) \, \widetilde{S}(G,Q',t) \\
& & - \sum_{Q\leq G' < P} g([G',P]^*,t) \sum_{\substack{Q'\leq E <
G'\\ Q \not\leq E}} g([E,G'],[E,E\vee Q],t)\,
\widetilde{S}(E,Q',t)\\
&=& \widetilde{S}(P,Q',t)+ \sum_{\substack{Q'\leq G < P\\ Q
\not\leq G}} g([G,P]^*,t) \, \widetilde{S}(G,Q',t) \\
& & - \sum_{\substack{Q'\leq E < P\\ Q \not\leq E}}
\widetilde{S}(E,Q',t) \sum_{E\vee Q \leq G' < P} g([E,G'],[E,E\vee
Q],t)\, g([G',P]^*,t)\\
&=& \widetilde{S}(P,Q',t)+ \sum_{\substack{Q'\leq F < P\\ Q
\not\leq F}}g([F,P],[F,F\vee Q],t)\, \widetilde{S}(F,Q',t),
\end{eqnarray*}
where the last step uses the formula for the relative
$g$-polynomial of 2.7. \hfill $\blacksquare$\\

\noindent \textbf{2.15. Corollary.} \textsl{If $Q'\leq Q$ are
faces of a lattice polytope $P$ then
\[ \widetilde{S}(P,Q',t) \leq \widetilde{S}(P,Q,t) \]
(i.e.\ the inequality holds coefficientwise). In particular,
$\widetilde{S}(P,Q,t)$ has nonnegative coefficients.}\\

\noindent \textit{Proof.} Note that the second statement follows
immediately since we can take $Q'= \emptyset$ and use that
$\widetilde{S}(P,t)$ has nonnegative coefficients. The first
statement is easily proved using induction on $\dim P$,
Proposition 2.14 and the nonnegativity of the coefficients of the
relative $g$-polynomials.
\hfill $\blacksquare$\\

\noindent \textbf{2.16.} For arbitrary order ideals $\mathcal{I}\varsubsetneq \mathcal{P}(P)$ of the face poset of a polytope $P$ we have the following recursion formula.\\

\noindent \textbf{Proposition.} \[\widetilde{S}(P,\mathcal{I},t) = h_P^*(t) - \sum_{\substack{\emptyset \leq F < P\\ F\notin \mathcal{I}}} \widetilde{S}(F,\mathcal{I}\cap \mathcal{P}(F) ,t)\, g([F,P]^*,t).  \] 

\noindent \textit{Proof.} The right hand side equals
\begin{eqnarray*}
& & h_P^*(t) - \sum_{\substack{\emptyset \leq F < P\\ F\notin \mathcal{I}}} \
\sum_{\substack{\emptyset \leq F' \leq F \\ F' \notin \mathcal{I} }} (-1)^{\dim F - \dim F' } h_{F'}^*(t)\, g([F',F],t)\, g([F,P]^*,t)\\
&=& h_P^*(t) - \sum_{\substack{\emptyset \leq F' < P\\ F' \notin \mathcal{I} }} h_{F'}^*(t)
\sum_{F' \leq F < P } (-1)^{\dim F - \dim F' } g([F',F],t)\, g([F,P]^*,t)\\
&=& h_P^*(t) + \sum_{\substack{\emptyset \leq F' < P\\ F' \notin \mathcal{I} }} (-1)^{\dim P - \dim F'} h_{F'}^*(t)\, g([F',P],t)\\ 
&=&  \widetilde{S}(P,\mathcal{I},t),
\end{eqnarray*}
where we used Stanley's convolution property for the $g$-polynomial. \hfill $\blacksquare$




\section{Stringy Hodge numbers}
\noindent \textbf{3.1.} In this section we review Batyrev's definition of stringy Hodge numbers, generalising
usual Hodge numbers of smooth projective varieties.\\

\noindent \textbf{3.2.} Let $X$ be a reduced but not necessarily
irreducible complex algebraic variety of dimension $d$. The
\textit{Hodge-Deligne polynomial} of $X$ is defined as
\[ H(X;u,v) := \sum_{i=0}^{2d} \sum_{p,q=0}^d
(-1)^{i} h^{p,q}(H_c^i(X,\mathbb{C}))\, u^pv^q  \quad \in
\mathbb{Z}[u,v], \] where $h^{p,q}(H_c^i(X,\mathbb{C}))$ denotes
the dimension of the $H^{p,q}$-component of the natural mixed
Hodge structure on $H_c^i(X,\mathbb{C})$. For a nonreduced variety $X$ we put $H(X;u,v):=H(X_{red};u,v)$. For a smooth projective
variety $X$, the coefficient of $u^pv^q$ in $H(X;u,v)$ is modulo
the factor $(-1)^{p+q}$ simply the Hodge number $h^{p,q}(X)$. The
Hodge-Deligne polynomial is a generalised Euler characteristic:
\begin{itemize}
\item[(1)] if $Y$ is a Zariski-closed subvariety of $X$, then $H(X;u,v) = H(Y;u,v) + H(X\setminus Y;u,v)$,
\item[(2)] for a product $X\times X'$ one has $H(X\times X';u,v) = H(X;u,v)\cdot H(X';u,v)$.
\end{itemize}

\vspace{0.3cm}

\noindent \textbf{3.3.} Let $Y$ from now on be a normal
irreducible variety. It is called $\mathbb{Q}$-\textit{Gorenstein}
if a multiple $rK_Y$ of the canonical class $K_Y$ is Cartier
($r\in \mathbb{Z}_{>0}$) and \textit{Gorenstein} if $K_Y$ itself
is Cartier. For example a normal hypersurface in a smooth variety
is Gorenstein. Let $f:X\to Y$ be a log resolution of $Y$. This
means that $f$ is a proper birational morphism from a smooth
variety $X$, such that the exceptional locus $D$ of $f$ is a
divisor with smooth components and normal crossings. Denote the
irreducible components of $D$ by $D_i$, where $i$ lives in a
finite index set $I$. For a $\mathbb{Q}$-Gorenstein $Y$ we have a linear equivalence
\[ rK_X \equiv f^*(rK_Y) + \sum_{i\in I} b_iD_i, \]
with $b_i\in \mathbb{Z}$ uniquely determined. One divides this formally by $r$ and calls the rational
number $a_i:= b_i/r$ the \textit{discrepancy coefficient} of $D_i$. The variety $Y$ is called
\textit{log terminal, canonical} or \textit{terminal} if all $a_i > -1, \geq 0$ or $>0$ respectively.
These definitions do not depend on the chosen log resolution and intuitively speaking these classes of
singularities are rather `mild'. If $Y$ is canonical but not terminal it is called \textit{strictly canonical}.\\

\noindent \textbf{3.4.} Let $Y$ be log terminal. Choose a log resolution $f:X\to Y$ and use the same
notations as above. For a subset $J$ of $I$ we set $D_J = \cap_{j\in J} D_j$ (so $D_{\emptyset}=X$)
and $D_J^{\circ} = D_J \setminus \cup_{i\in I\setminus J} D_i$. The varieties $D_J^{\circ}$ give a
natural stratification of $X$.\\

\noindent \textbf{Definition} \cite[Def.\ 3.1]{Batyrev}\textbf{.}
The \textit{stringy $E$-function} of $Y$ is defined as
\[ E_{st}(Y;u,v) := \sum_{J\subset I} H(D_J^{\circ};u,v) \prod_{j\in J} \frac{uv-1}{(uv)^{a_j+1}-1}, \]
where $a_j$ is the discrepancy coefficient of $D_j$ and where the product over $j$ has to be interpreted as 1 if $J=\emptyset$.\\

\noindent Batyrev used motivic integration to prove that this formula does not depend on the chosen resolution
\cite[Thm.\ 3.4]{Batyrev}.\\

\noindent \textbf{3.5.} The following remarks are in order:
\begin{itemize}
\item[(1)] If $Y$ is Gorenstein, then $E_{st}(Y;u,v)$ is a rational function. It lives then
in $\mathbb{Z}[[u,v]] \cap \mathbb{Q}(u,v)$.
\item[(2)] If $Y$ is smooth, then $E_{st}(Y;u,v) = H(Y;u,v)$. More generally, if $Y$ has a crepant resolution (i.e.\ a
log resolution $f:X\to Y$ such that all discrepancy coefficients are 0), then $E_{st}(Y;u,v)=H(X;u,v)$.
\item[(3)] We can choose the log resolution $f:X\to Y$ such that it is an isomorphism when restricted to the inverse
image of the nonsingular part of $Y$. In particular, using such a log resolution, we see that if $Y$ has only an
isolated singularity at a point $y$, then we can write
\[ E_{st}(Y;u,v) = H(Y\setminus \{y\};u,v) + \sum_{\emptyset\neq J\subset I} H(D_J^{\circ};u,v) \prod_{j\in J}
\frac{uv-1}{(uv)^{a_j+1}-1}. \]
We call $E_{st}(Y;u,v) - H(Y\setminus \{y\};u,v)$ the \textit{local contribution} of the isolated singularity and
denote it by $E_{st,y}(Y;u,v)$.
\end{itemize}

\vspace{0.3cm}

\noindent \textbf{3.6.} For a projective variety $Y$ of dimension
$d$ Batyrev \cite[Thm.\ 3.7]{Batyrev} proved that
\begin{itemize}
\item[(1)] $E_{st}(Y;u,v)= (uv)^d E_{st}(Y;u^{-1},v^{-1}),$
\item[(2)] $E_{st}(Y;0,0) = 1$.
\end{itemize}
Note that this generalises the relations $h^{p,q}(Y)=h^{q,p}(Y)=h^{d-p,d-q}(Y)=h^{d-q,d-p}(Y)$ and $h^{0,0}(Y)=1$
for a smooth projective $Y$.\\

\noindent \textbf{3.7.} Let $Y$ now be a projective variety of
dimension $d$ with Gorenstein canonical singularities. Assume that
$E_{st}(Y;u,v)$ is a polynomial $\sum_{p,q} a_{p,q}u^pv^q$.
Batyrev then defines the \textit{stringy Hodge numbers} of $Y$ as
$h^{p,q}_{st}(Y) = (-1)^{p+q}a_{p,q}$ \cite[Def.\ 3.8]{Batyrev}.
By 3.5 (2) and 3.6 one has
\begin{itemize}
\item[(1)] stringy Hodge numbers $h_{st}^{p,q}(Y)$ can only be nonzero for $0\leq p\leq d$ and $0\leq q \leq d$,
\item[(2)] for smooth projective varieties stringy Hodge numbers are equal to usual Hodge numbers,
\item[(3)] $h_{st}^{p,q}(Y)=h_{st}^{q,p}(Y)=h_{st}^{d-p,d-q}(Y)=h_{st}^{d-q,d-p}(Y)$ and $h_{st}^{0,0}(Y)=1$.
\end{itemize}
The following intriguing question is however still open.\\

\noindent \textbf{3.8. Conjecture} \cite[Conj.\ 3.10]{Batyrev}\textbf{.} \textsl{Stringy Hodge
numbers are nonnegative.}\\

\noindent We remark that it is not clear when to expect a polynomial stringy $E$-function. It is true for Gorenstein
toric varieties \cite[Prop.\ 4.4]{Batyrev} and the stringy Hodge numbers (or in that case better stringy Betti
numbers) are nonnegative for complete Gorenstein toric varieties \cite[Thm.\ 1.2]{MustataPayne}. For varieties with
Gorenstein quotient singularities the stringy $E$-function is also polynomial. Yasuda showed that the stringy Hodge
numbers for such complete varieties coincide with the orbifold cohomology Hodge numbers \cite[Rem.\ 1.4 (2)]{Yasuda}.
In the next section we describe a natural class of isolated strictly canonical nondegenerate hypersurface singularities
that also have a polynomial stringy $E$-function. In Section 5 we prove that Batyrev's conjecture holds for complete
varieties with such singularities, under the additional hypothesis of Theorem 1.4.\\

\noindent Batyrev's conjecture is easy for surfaces. Indeed,
canonical surface singularities are classified: it are precisely
the so-called $A$-$D$-$E$ singularities. One knows that these
singularities admit a crepant resolution and hence Batyrev's
conjecture for surfaces follows from 3.5 (2). In higher dimension
there is the following theorem
\cite[Thm.\ 3.1 and Cor.\ 3.4]{SchepersVeys}.\\

\noindent \textbf{3.9. Theorem.} \begin{itemize}
\item[(1)] \textsl{For threefolds, Batyrev's conjecture is true in full
generality.}
\item[(2)] \textsl{Let $Y$ be a projective variety of dimension $d\geq 4$ with at
most isolated Gorenstein canonical singularities and with
polynomial stringy $E$-function. Assume that $Y$ has a log
resolution $f:X\to Y$ such that all discrepancy coefficients of
irreducible exceptional components are $> \lfloor \frac{d-4}{2}
\rfloor$. Then the stringy Hodge numbers of $Y$ are nonnegative.}
\end{itemize}

\section{Nondegenerate singularities}
\noindent \textbf{4.1.} In this section we recall the definition of nondegenerate hypersurface singularities and we
explain how to compute their stringy $E$-function. We describe a natural class of isolated strictly canonical
nondegenerate singularities that give rise to a polynomial stringy $E$-function. We conclude by giving a concrete formula for this contribution to the stringy $E$-function.\\

\noindent \textbf{4.2.} Let $f\in \mathbb{C}[x_1,\ldots,x_n]$ be a
polynomial with $f(\mathbf{0})=0$. We denote the hypersurface
$\{f=0\}\subset \mathbb{C}^n$ by $X_f$. We write
$f=\sum_{\mathbf{m}\in (\mathbb{Z}_{\geq 0})^n} a_{\mathbf{m}}
\mathbf{x}^{\mathbf{m}}$ where $\mathbf{x}^{\mathbf{m}} =
x_1^{m_1}\cdots x_n^{m_n}$. The \textit{Newton polyhedron}
$\Gamma(f)$ of $f$ is the convex hull in $\mathbb{R}^n$ of
\[ \bigcup_{\substack{\mathbf{m}\in (\mathbb{Z}_{\geq 0})^n   \\ a_{\mathbf{m}}\neq 0 }} \mathbf{m} +
(\mathbb{R}_{\geq 0})^n.   \]
A \textit{face} of $\Gamma(f)$ is defined as any nonempty intersection of $\Gamma(f)$ with a hyperplane $H$ such that
$\Gamma(f)$ is completely contained in one of the two closed halfspaces determined by $H$. This is similar to the
definition of a face of a polytope, but now we do not consider the empty set as a face and hence $\Gamma(f)$ is the
only improper face of $\Gamma(f)$. For a face $\tau$ of $\Gamma(f)$ we write $f_{\tau}$ for the polynomial
$\sum_{\mathbf{m}\in \tau \cap (\mathbb{Z}_{\geq 0})^n} a_{\mathbf{m}} \mathbf{x}^{\mathbf{m}}$. One calls $f$
\textit{nondegenerate with respect to its Newton polyhedron} if the equation $f_{\tau}=0$ defines a smooth subvariety
of $(\mathbb{C}^*)^n$ for every \textsl{compact} face $\tau$ of $\Gamma(f)$.\\

\noindent \textbf{4.3.} From the Newton polyhedron of a polynomial
$f$ one gets a partition of $(\mathbb{R}_{\geq 0})^n$ into cones
(where $(\mathbb{R}_{\geq 0})^n$ should be considered as the first
orthant of the space dual to the surrounding space of
$\Gamma(f)$). This goes as follows. For a vector $\mathbf{v}\in
(\mathbb{R}_{\geq 0})^n$ set $m_f(\mathbf{v}) =
\inf_{\mathbf{w}\in \Gamma(f)}\{\mathbf{v}\cdot \mathbf{w}\}$,
where $\cdot$ is the standard inner product. In fact this infimum
is attained and hence it is a minimum. The \textit{first meet locus}
$F(\mathbf{v})$ of $\mathbf{v}$ is defined as
\[ F(\mathbf{v}):=\{ \mathbf{w} \in \Gamma(f)\,|\,  \mathbf{v}\cdot \mathbf{w} =m_f(\mathbf{v}) \}.\]
This is a face of $\Gamma(f)$ and it is a compact face if and only if $\mathbf{v}\in (\mathbb{R}_{>0})^n$. For a face
$\tau$ of $\Gamma(f)$ we can then define the \textit{cone $\delta_{\tau}$ associated to} $\tau$ by
\[ \delta_{\tau}:= \{ \mathbf{v} \in (\mathbb{R}_{\geq 0})^n\,|\, F(\mathbf{v})=\tau  \}.\]
These cones form a partition of $(\mathbb{R}_{\geq 0})^n$ and their closures are pointed rational convex polyhedral
cones with vertex at the origin, forming a fan $\Delta_f$. Following \cite[\S 5]{Stepanov} we call this fan $\Delta_f$
the \textit{first Varchenko subdivision}. For a nondegenerate $f$ such that the origin is an isolated singularity
of $X_f$, this construction gives the first step in a toric resolution of $(X_f,\mathbf{0})$ \cite[\S 9, 10]{Varchenko}. More
precisely, $\Delta_f$ can be subdivided to a fan $\Delta'$ consisting of unimodular cones (i.e.\ simplicial cones that
can be generated by a part of a $\mathbb{Z}$-basis of $\mathbb{Z}^n$). Then the proper birational toric map from the
toric variety $X(\Delta')$, associated to $\Delta'$, to $\mathbb{C}^n$ is an embedded resolution of singularities of
$(X_f,\mathbf{0})$.\\

\noindent \textbf{4.4.} Let $g=\sum_{\mathbf{m}\in \mathbb{Z}^n} b_{\mathbf{m}}
\mathbf{x}^{\mathbf{m}}\in \mathbb{C}[x_1,\ldots,x_n,x_1^{-1},\ldots,x_n^{-1}]$ be a Laurent polynomial. The
\textit{Newton polytope} $P$ of $g$ is the convex hull in $\mathbb{R}^n$ of those $\mathbf{m}\in \mathbb{Z}^n$ with
$b_{\mathbf{m}}\neq 0$. It is a lattice polytope. For a nonempty face $F$ of $P$ we write $g_F$ for the Laurent
polynomial $\sum_{\mathbf{m}\in F\cap \mathbb{Z}^n} b_{\mathbf{m}} \mathbf{x}^{\mathbf{m}}$. One calls $g$
\textit{nondegenerate with respect to its Newton polytope} if the equation $g_F=0$ defines a smooth subvariety in
$(\mathbb{C}^*)^n$ for every nonempty face $F$ of $P$. \\

\noindent \textbf{4.5.} Let $P$ be a lattice polytope in $\mathbb{R}^n$ of maximal dimension. Let $g$ be a Laurent
polynomial with $P$ as Newton polytope and assume that $g$ is nondegenerate with respect to $P$. Batyrev and Borisov
derived a formula for the Hodge-Deligne polynomial of the hypersurface $Y_g:=\{g=0\}\subset (\mathbb{C}^*)^n$
\cite[Thm.\ 3.18]{BatyrevBorisov}. In the proof of Proposition 5.5 of \cite{BorisovMavlyutov} this formula is rewritten
as follows.\\

\noindent \textbf{Theorem.} \textsl{Using the notations of Section
2, $H(Y_g;u,v) $ equals}
\[ \frac{1}{uv}\biggl( (uv-1)^{\dim P} +(-1)^{\dim P +1} \sum_{\emptyset\leq F\leq P} u^{\dim F+1}\,
\widetilde{S}(F,u^{-1}v)\,g([F,P]^*,uv) \biggr).  \]

\noindent \textbf{4.6.} Now let $f\in \mathbb{C}[x_1,\ldots,x_n]$ be an irreducible polynomial with
$f(\mathbf{0})=0$ that is nondegenerate with respect to its Newton polyhedron. Assume that the hypersurface
$X_f$ has only an isolated singularity at $\mathbf{0}$. For a vector $\mathbf{v}=(v_1,\ldots,v_n)\in \mathbb{R}^n$
set $\sigma(\mathbf{v}):=v_1+\cdots+v_n$. \\

\noindent \textbf{Proposition.}
\begin{itemize}
\item[(1)] $X_f$ is canonical $\Leftrightarrow$ for all primitive vectors $\mathbf{v}\in (\mathbb{Z}_{\geq 0})^n$
we have $\sigma(\mathbf{v})-m_f(\mathbf{v})\geq 1$ (primitive means that $\gcd(v_1,\ldots,v_n)=1$).
\item[(2)] $X_f$ is terminal $\Leftrightarrow$ for all primitive vectors $\mathbf{v}\in (\mathbb{Z}_{\geq 0})^n$
different from the standard basis vectors $\mathbf{e}_i$ ($i=1,\ldots,n$) we have $\sigma(\mathbf{v})-
m_f(\mathbf{v})> 1$.
\end{itemize}

\noindent \textit{Proof.} This is essentially Theorem 4.6 of
\cite{YPG}, but two remarks are in order. Firstly, the phrase
`different from the standard basis vectors' is missing in the
formulation of that theorem but it should be there. Secondly, we
have to explain that the above conditions for being canonical or
terminal are not only necessary but also sufficient in the
nondegenerate case. Let $\Delta'$ be a fan as in 4.3. Let $\delta$
be a cone of maximal dimension of $\Delta'$ generated by integer
vectors $\delta^1=(\delta^1_1,\ldots , \delta^1_n),\ldots,\delta^n=(\delta^n_1,\ldots , \delta^n_n)$ that form a $\mathbb{Z}$-basis
for $\mathbb{Z}^n$. The affine toric variety $X(\delta)$
associated to $\delta$ is isomorphic to $\mathbb{C}^n$ and
$X(\Delta')$ is covered by the open sets of the form $X(\delta)$.
The proper birational toric map $h:X(\Delta') \to \mathbb{C}^n$ is
locally given by
\[ h: X(\delta)\cong \mathbb{C}^n \to \mathbb{C}^n: (y_1,\ldots,y_n) \mapsto (x_1 = \prod_i y_i^{\delta^i_1},
\ldots,x_n = \prod_i y_i^{\delta^i_n})\] and the total inverse
image of $X_f$ on $X(\delta)$ is given by an equation
\[ y_1^{m_f(\delta^1)}\cdots y_n^{m_f(\delta^n)} f_{\delta}(y_1,\ldots,y_n) =0,   \]
with $f_{\delta}(0,\ldots,0)\neq 0$ \cite[\S 10]{Varchenko}. Since $f$ is nondegenerate, $h$ gives an embedded resolution
of singularities of $X_f$. Let $X_f'$ be the proper transform of $X_f$ under $h$. The discrepancy coefficient of an
exceptional component of the induced log resolution $h:X_f'\to X_f$ can be computed from the embedded resolution using
the adjunction formula (for details see for instance the proof of Proposition 2.3 of \cite{SchepersVeys2}). It equals
$\sigma(\delta^i)-m_f(\delta^i)-1$ for an exceptional component whose intersection with $X(\delta)$ is nonempty and
lies in $\{y_i=0\}$. \hfill $\blacksquare$\\

\noindent \textbf{4.7.} Assume now that $X_f$ is canonical and has only an isolated singularity at the origin. The
local contribution $E_{st,\mathbf{0}}(X_f;u,v)$ of the singularity to the stringy $E$-function of $X_f$ can be
computed by the following result (\cite[Cor.\ 3.2]{SchepersVeys2}, essentially work of Denef and Hoornaert
\cite{DenefHoornaert}).\\

\noindent \textbf{Proposition.} \textsl{
\[ E_{st,\mathbf{0}}(X_f;u,v) = \sum_{\substack{\text{compact faces}\\ \tau\text{ of }\Gamma(f)}} H(N_{\tau};u,v)\,
T_{\delta_{\tau}}(f,uv), \]
where $N_{\tau}$ is the subvariety of $(\mathbb{C}^*)^n$ given by $\{f_{\tau}=0\}$ and where $T_{\delta_{\tau}}(f,t)$
is the power series $\sum_{\mathbf{v} \in (\mathbb{Z}_{>0})^n
\cap \delta_{\tau}}
t^{m_f(\mathbf{v})-\sigma(\mathbf{v})}$ (one can show that this power series belongs to $\mathbb{Q}(t)$).}\\

\noindent \textbf{4.8.} We are ready to define the singularities that we study during the rest of this paper.\\

\noindent \textbf{Definition.} Let $f\in \mathbb{C}[x_1,\ldots,x_n]$ be nondegenerate with respect to its Newton
 polyhedron and let $X_f$ have an isolated singularity at $\mathbf{0}$. We call the first Varchenko subdivision
 $\Delta_f$ \textit{crepant} if all primitive integer generators $\delta^i$ of 1-dimensional cones of $\Delta_f$ satisfy
 $\sigma(\delta^i)-m_f(\delta^i)=1$.\\

\noindent \textbf{4.9. Remark.}
\begin{itemize}
\item[(1)] Let $\Delta'$ be a subdivision of $\Delta_f$ in unimodular cones. As explained above, this subdivision
gives an embedded resolution $X(\Delta')\to \mathbb{C}^n$ of $X_f$, inducing a log resolution. If $\Delta_f$ is
crepant then the exceptional components of this log resolution coming from 1-dimensional cones of $\Delta_f$ all
have discrepancy coefficient 0, so this explains the name.
\item[(2)] Note that $\sigma-m_f$ is a linear function when restricted to a cone $\delta$ of $\Delta_f$. So if
$\delta$ is of maximal dimension then $(\sigma-m_f)|_{\delta}$ gives rise to an element $n_{\delta}$ of the dual
Hom$(\mathbb{Z}^n,\mathbb{Z})$ of $\mathbb{Z}^n$. If $\Delta_f$ is crepant, then using these $n_{\delta}$ we see
that every cone $\delta$ of maximal dimension is a \textit{Gorenstein cone} (see \cite[Def.\ 1.8]{BatyrevNill}).
\item[(3)] Proposition 4.6 shows that these singularities are strictly canonical: any primitive integer vector
$\mathbf{v}$ is a linear combination of primitive generators of a cone of $\Delta_f$ with nonnegative rational
coefficients and hence $\sigma(\mathbf{v})-m_f(\mathbf{v})\geq 1$.
\end{itemize}

\noindent \textbf{4.10. Example.} In the following table we
investigate the condition of having a crepant first Varchenko
subdivision for the standard equations of the canonical or
$A$-$D$-$E$ surface singularities. In the second column we give
the primitive generators of the 1-dimensional cones of $\Delta_f$
different from the standard basis vectors.

\setlongtables
\begin{longtable}{c|c|c}
& & \\ [-0.3cm] Singularity & Primitive generators  & $\Delta_f$
crepant\,?  \\ [0.2cm] \hline & & \\ [-0.3cm] $\begin{array}{c}
A_n: x^{n+1}+y^2+z^2=0\\ (n\geq 1) \end{array}$ &
$\begin{array}{l} n \text{ even}: (2,n+1,n+1)\\ n \text{ odd}:
(1,\frac{n+1}{2},\frac{n+1}{2}) \end{array}$ & $\begin{array}{c}
\text{no}\\ \text{yes} \end{array}$ \\ [0.4cm]
 \hline  & & \\ [-0.3cm]
$\begin{array}{c} D_n:x^{n-1}+xy^2 + z^2=0\\ (n\geq 4) \end{array}$ & $(2,n-2,n-1),(2,0,1)$ & yes \\ [0.3cm]
\hline & & \\ [-0.3cm]
$E_6: x^4+y^3+z^2=0$ & $(3,4,6)$ & yes  \\ [0.2cm]
\hline  & &   \\ [-0.3cm]
$E_7: x^3+xy^3+z^2=0$ & $(6,4,9),(2,0,1)$ & yes\\ [0.2cm]
 \hline & & \\ [-0.3cm]
$E_8: x^5 + y^3 + z^2=0$ & $(6,10,15)$ & yes \\ [0.1cm]
\end{longtable}

\noindent \textbf{4.11. Proposition.} \textsl{Let $f\in \mathbb{C}[x_1,\ldots,x_n]$ be nondegenerate with respect
to its Newton polyhedron and let $X_f$ have an isolated singularity at $\mathbf{0}$. If the first Varchenko
subdivision is crepant, then the local contribution $E_{st,\mathbf{0}}(X_f;u,v)$ is a polynomial.}\\

\noindent \textit{Proof.} We use Proposition 4.7. First we rewrite
$T_{\delta_{\tau}}(f,uv)$. Let $\tau$ be a compact face of
$\Gamma(f)$, $\delta_{\tau}$ the associated cone and
$\overline{\delta_{\tau}}$ its closure. Let $P_{\tau}$ be the
convex hull in $\mathbb{R}^n$ of the primitive integer generators
of the extreme rays of $\overline{\delta_{\tau}}$. Let $\tau'$ be
a vertex of $\tau$. Then $\overline{\delta_{\tau'}}$ is a cone of
maximal dimension of $\Delta_f$. Let $n_{\delta_{\tau'}}$ be as in
Remark 4.9 (2). Then $P_{\tau'}$ is the so-called \textit{support
polytope} of the Gorenstein cone $\overline{\delta_{\tau'}}$
(i.e.\ all points where $n_{\delta_{\tau'}}$ takes value 1), and
$P_{\tau} = P_{\tau'} \cap \overline{\delta_{\tau}}$. We look at
$n_{\delta_{\tau'}}$ as a degree function on $P_{\tau'}$ and on
$P_{\tau}$ (it obviously does not depend on the choice of the
vertex). Then we have \[ T_{\delta_{\tau}}(f,uv)=\sum_{\mathbf{v}
\in (\mathbb{Z}_{>0})^n \cap \delta
 _{\tau}} (uv)^{-n_{\delta_{\tau'}}(\mathbf{v})}. \]
By Stanley's reciprocity law (formulated in \cite[Thm.\ 4.6.14]{Stanley97} for solutions of linear homogeneous
diophantine equations, but allowing inequalities is no problem by the remark on p.222 of \cite{Stanley97}) this equals
\[ (-1)^{\dim \delta_{\tau}}  \sum_{\mathbf{v} \in (\mathbb{Z}_{\geq 0})^n \cap \overline{\delta_{\tau}}}
(uv)^{n_{\delta_{\tau'}}(\mathbf{v})}.\]
This can be rewritten as
\[(-1)^{\dim \delta_{\tau}} \sum_{i\in \mathbb{Z}_{\geq 0}}  | (iP_{\tau})\cap \mathbb{Z}^n | \, (uv)^i =
\frac{h^*_{P_{\tau}}(uv) }{(uv-1)^{\dim \delta_{\tau}}}\]
by definition of the $h^*$-polynomial of the lattice polytope $P_{\tau}$ of dimension $\dim \delta_{\tau}-1$.
Using Proposition 4.7 it suffices now to show that $\frac{H(N_{\tau};u,v)}{(uv-1)^{\dim \delta_{\tau}}}$ is a
polynomial for each compact face $\tau$. First we divide the equation $f_{\tau}$ by one of the monomials appearing
in it (this corresponds to moving one of the vertices of $\tau$ to the origin). We get a Laurent polynomial
$\widetilde{f_{\tau}}$ and $N_{\tau} \cong (\{\widetilde{f_{\tau}}=0\}\subset (\mathbb{C}^*)^n)$. Then we use a
coordinate change
\[     y_j = \prod_{i=1}^n x_i^{t_{i,j}},\text{ where } j=1,\ldots,n \text{ and } T=(t_{i,j})\in GL_n(\mathbb{Z}),\]
on $(\mathbb{C}^*)^n$ as in Lemma 5.8 of \cite{DenefHoornaert} to write
\[  \widetilde{f_{\tau}}(x_1,\ldots,x_n) = h_{\tau}(y_1,\ldots,y_{\dim \tau}) \]
for a (nonunique) Laurent polynomial $h_{\tau}$. Then
\[N_{\tau} \cong (\mathbb{C^*})^{n-\dim \tau} \times (\{h_{\tau} = 0 \} \subset (\mathbb{C}^*)^{\dim \tau}).\]
Since $H(\mathbb{C}^*;u,v)=uv-1$ and $n-\dim \tau = \dim \delta_{\tau}$ we conclude that
$\frac{H(N_{\tau};u,v)}{(uv-1)^{\dim \delta_{\tau}}}$ is a polynomial.\hfill $\blacksquare$\\

\noindent \textbf{4.12. Example.} The previous proposition is in general not true for isolated strictly canonical nondegenerate singularities that do not have a crepant first Varchenko subdivision. Consider the polynomial $f = x_1^5 + x_2^3 + x_3^3 + x_4^3$. Proposition 4.3 from \cite{Reid} shows that the singularity $(X_f,\mathbf{0})$ is canonical and Proposition 4.6 applied with $\mathbf{v} = (1,1,1,1)$ shows that it is strictly canonical. Using the combinatorial procedure of \cite[Section 4]{SchepersVeys2} one finds the expression
\[ \frac{(uv)^4-u^4v^3-u^3v^4+3(uv)^3-2u^3v^2-2u^2v^3+4(uv)^2-u^2v-uv^2+2uv+1}{(uv)^2+uv+1}\]
for $E_{st,\mathbf{0}}(X_f;u,v)$.\\

\noindent \textbf{4.13.} In the proof of Proposition 4.11 we associated a lattice polytope $P_{\tau}$ to a compact face $\tau$ of the Newton polyhedron
$\Gamma(f)$ of $f$. Now we define the lattice polytope $P_{\emptyset}$ as the
convex hull of all the polytopes $P_{\tau}$ and the origin in
$(\mathbb{R}_{\geq 0})^n$. So $P_{\emptyset}$ is a kind of
fundamental domain of the crepant first Varchenko subdivision
$\Delta_f$. Denote the set of compact faces of $\Gamma(f)$ by $\mathcal{P}_f$. We get an inclusion-reversing bijective correspondence between $\mathcal{P}_f\cup \{\emptyset\}$ and the faces of $P_{\emptyset}$ that are not contained in a coordinate hyperplane. Using these notations and the notations of Section 2 we can summarise the results so far in the following theorem.\\

\noindent \textbf{4.14. Theorem.} \textsl{Let $f\in \mathbb{C}[x_1,\ldots,x_n]$ be nondegenerate with respect
to its Newton polyhedron and let $X_f$ have an isolated singularity at $\mathbf{0}$. If the first Varchenko
subdivision is crepant, then the local contribution $E_{st,\mathbf{0}}(X_f;u,v)$ is given by 
\[ \frac{1}{uv} \sum_{\mu \in \mathcal{P}_f\cup \{\emptyset\}} (-u)^{\dim \mu +1}\,
\widetilde{S}(\mu,u^{-1}v)\,
\widetilde{S}(P_{\mu},star_{\partial P_{\emptyset}}(\mathbf{0})\cap \mathcal{P}(P_{\mu}),uv) .\]}

\noindent \textit{Proof.} In the proof of Proposition 4.11 we showed that
\[N_{\tau} \cong (\mathbb{C^*})^{n-\dim \tau} \times (\{h_{\tau} = 0 \} \subset (\mathbb{C}^*)^{\dim \tau})\]
for a Laurent polynomial $h_{\tau}$ in variables $y_1,\ldots,
y_{\dim \tau}$. Since $f$ is nondegenerate with respect
to its Newton polyhedron, $h_{\tau}$ is nondegenerate with respect
to its Newton polytope in $\mathbb{Z}^{\dim \tau}$. If we consider
$\tau$ as a lattice polytope as well, then this Newton polytope is
by construction isomorphic to $\tau$. Combining Proposition 4.7,
the proof of Proposition 4.11 and Theorem 4.5 we get
\[ E_{st,\mathbf{0}}(X_f;u,v) = \sum_{\tau \in \mathcal{P}_f} \frac{h^*_{P_{\tau}}(uv)}{uv}
\biggl( (uv-1)^{\dim \tau} +(-1)^{\dim \tau +1} A(\tau;u,v) \biggr),\]
where
\[ A(\tau;u,v) = \sum_{\emptyset\leq \mu\leq \tau} u^{\dim \mu +1} \widetilde{S}(\mu,u^{-1}v)\, g([\mu,\tau]^*,uv).\]
We split $E_{st,\mathbf{0}}(X_f;u,v)$ as $A_1(uv)+A_2(u,v)$ with
\begin{eqnarray*}
A_1(uv) & := &  \sum_{\tau \in \mathcal{P}_f} \frac{h^*_{P_{\tau}}(uv)}{uv}\biggl( (uv-1)^{\dim \tau} +
(-1)^{\dim \tau +1} g([\emptyset, \tau]^*,uv) \biggr) , \\
A_2(u,v) & := & \sum_{\tau \in \mathcal{P}_f} (-1)^{\dim \tau
+1}\, \frac{h^*_{P_{\tau}}(uv)}{uv}\,  \overline{A}(\tau;u,v),
\end{eqnarray*}
where
\[ \overline{A}(\tau;u,v) = \sum_{\emptyset < \mu\leq \tau} u^{\dim \mu +1}
\widetilde{S}(\mu,u^{-1}v)\, g([\mu,\tau]^*,uv).\] 
We remark that
$A_1(uv)$ and $A_2(u,v)$ are both polynomials: the constant
coefficient of a $g$-polynomial is always 1, so it follows
immediately that $A_1(uv)$ is a polynomial. And $u^{\dim \mu +1}
\widetilde{S}(\mu,u^{-1}v)$ is a homogeneous polynomial in $u,v$
of degree $\dim \mu +1$ without terms of the form $cu^{\dim \mu
+1}$ or $cv^{\dim \mu +1}$ by
2.12. This shows that $A_2(u,v)$ is a polynomial.\\

\noindent Note that $P_{\emptyset}$ with chosen vertex
$\mathbf{0}$ satisfies the conditions of Lemma 2.10. Using that
lemma and Definition 2.13 we rewrite $A_1(uv)$ as
\[  \frac{1}{uv}\bigl(h_{P_{\emptyset}}^*(uv) + \widetilde{S}(P_{\emptyset},
star_{\partial P_{\emptyset}}(\mathbf{0}),uv)-h_{P_{\emptyset}}^*(uv) \bigr) = \frac{\widetilde{S}(P_{\emptyset},star_{\partial P_{\emptyset}}(\mathbf{0}),uv)}{uv}.  \tag{3}            \]
By exchanging the sums we find that $A_2(u,v)$ equals
\[\frac{1}{uv} \sum_{\mu \in \mathcal{P}_f} (-u)^{\dim \mu +1} \widetilde{S}(\mu,u^{-1}v)
 \sum_{\substack{\tau \in \mathcal{P}_f\\ \tau \geq \mu}} (-1)^{\dim\tau -\dim \mu}\, h_{P_{\tau}}^*(uv)\,
  g([\mu,\tau]^*,uv) \]
\[= \frac{1}{uv} \sum_{\mu \in \mathcal{P}_f} (-u)^{\dim \mu +1} \widetilde{S}(\mu,u^{-1}v)\,
 \widetilde{S}(P_{\mu},star_{\partial P_{\emptyset}}(\mathbf{0})\cap \mathcal{P}(P_{\mu}),uv). \tag{4}\]
Adding formulae (3) and (4) ends the proof of the theorem.\hfill $\blacksquare$

\section{Nonnegativity of stringy Hodge numbers}
\noindent \textbf{5.1.} In the previous section we defined a class of isolated nondegenerate singularities that
give a polynomial contribution to the stringy $E$-function. Now we will study Batyrev's conjecture about the nonnegativity of the stringy Hodge numbers for these singularities. First we formulate a lemma that allows us to draw (global) conclusions
about the stringy Hodge numbers from the local contributions of the singularities.\\

\noindent \textbf{5.2. Lemma.} \textsl{Let $Y$ be a complete variety of dimension $d$ with at most isolated
Gorenstein canonical singularities and with a polynomial stringy $E$-function. Assume that the local contribution
of the singularities to the stringy $E$-function is $\sum_{i,j} c_{i,j} u^iv^j$ with $(-1)^{i+j}c_{i,j} \geq 0$
for $i+j\geq d$. Then the stringy Hodge numbers of $Y$ are nonnegative.}\\

\noindent \textit{Proof.} This is a generalisation of Lemma 4.4 from \cite{SchepersVeys2}. The proof given there
proves in fact exactly this generalisation. \hfill $\blacksquare$ \\

\noindent \textbf{5.3. Theorem.} \textsl{Let $Y$ be a complete variety of dimension $d$ with isolated
singularities, such that each singularity is analytically isomorphic to a nondegenerate hypersurface singularity $(X_f,\mathbf{0})$ with crepant first
Varchenko subdivision. Assume in addition that each defining polynomial $f\in \mathbb{C}[x_1,\ldots,x_{d+1}]$ of an occurring nondegenerate singularity satisfies
\begin{enumerate}
\item $f$ is convenient (i.e.\ $f$ contains a nonzero term $a_ix_i^{b_i}$ for each variable $x_i$),
\item \textbf{or} $\Gamma(f)$ has a unique maximal compact face.
\end{enumerate}
Then the stringy Hodge numbers of $Y$ are nonnegative.}\\

\noindent \textbf{5.4. Remark.}
\begin{itemize}
\item[(1)] One should remark that these singularities in general do not allow a crepant
resolution (see 3.5 (2)), as Example 4.3 from \cite{Schepers} shows. Indeed, in that example the stringy Hodge
numbers do not satisfy the Hard Lefschetz property and hence they cannot be Hodge numbers of a smooth projective
variety.
\item[(2)] The second condition on $f$ holds for example if $f$ is weighted homogeneous.
\item[(3)] Theorem 5.3 and Proposition 4.11 form a strong generalisation of most of the results on Brieskorn
singularities from \cite[Section 4]{SchepersVeys2}.
\item[(4)] It is interesting to compare this result with Theorem 3.9. From the point of view of that theorem,
strictly canonical singularities are the worst untreated case. On
the other hand, strictly canonical singularities should provide
most of the examples of polynomial stringy $E$-functions since the
a priori denominator of the stringy $E$-function becomes worse if
the discrepancy coefficients are bigger. These are two good
reasons to study strictly canonical singularities in this context.
\end{itemize}

\noindent \textit{Proof of Theorem 5.3.} Let $y\in Y$ be a singular point of $Y$, analytically isomorphic to the
hypersurface singularity $(X_f,\mathbf{0})$ for a nondegenerate polynomial $f\in \mathbb{C}[x_1,\ldots,x_{d+1}]$ with crepant first
Varchenko subdivision. We will show that the local contribution $E_{st,y}(Y;u,v)=\sum_{i,j} c_{i,j}u^iv^j$ satisfies
$(-1)^{i+j}c_{i,j}\geq 0$ for $i+j\geq d$ under one of the extra conditions from the theorem on $f$. Then we can apply Lemma 5.2 to deduce the theorem.\\

\noindent We use the formula from Theorem 4.14. Note that $\frac{1}{uv} (-u)^{\dim \mu +1}
\widetilde{S}(\mu,u^{-1}v)$ is a homogeneous polynomial of degree $\dim\mu -1$ by 2.12 if $\mu\neq\emptyset$ (it is 0 if $\dim \mu =0$). The sign of its coefficients simply depends on the parity of the degree. If $\mu = \emptyset$ this expression equals $\frac{1}{uv}$, but recall from the proof of Theorem 4.14 that $\frac{\widetilde{S}(P_{\emptyset},star_{\partial P_{\emptyset}}(\mathbf{0}),uv)}{uv}$ is a polynomial. Hence it suffices to show that for all $\mu\in \mathcal{P}_f\cup \{\emptyset\}$ the expression 
\[ \tau_{>(d-\dim\mu)/2}\,\widetilde{S}(P_{\mu},\mathcal{I}_{\mu},t)  \] 
has nonnegative coefficients, where we wrote $\mathcal{I}_{\mu}$ for $star_{\partial P_{\emptyset}}(\mathbf{0})\cap \mathcal{P}(P_{\mu})$ and where $\tau_{>}$ denotes the truncation operator as in Definition 2.3.
\begin{enumerate}
\item First assume that $f$ is convenient. We will prove by descending induction on $\dim\mu$ that 
\[ \tau_{>(d-\dim\mu)/2}\,\widetilde{S}(P_{\mu},\mathcal{I}_{\mu},t)  = \tau_{>(d-\dim\mu)/2}\,\widetilde{S}(P_{\mu},t) . \tag{5} \] If $\dim \mu = d$ then $P_{\mu}$ is a vertex and the equality (5) is trivial. Assume now that $\dim \mu <d$. By Proposition 2.16 we have that $\tau_{>(d-\dim\mu)/2}\,\widetilde{S}(P_{\mu},\mathcal{I}_{\mu},t)$ equals
\[ \tau_{>(d-\dim\mu)/2}\, h_{P_{\mu}}^*(t) - \sum_{\substack{\nu \in \mathcal{P}_f  \\ \nu > \mu }}  \sum_{i=0}^{r_{\nu,\mu}} \tau_{=i}\,g([\mu,\nu],t) \cdot\tau_{> (d-\dim\mu)/2 -i}\, \widetilde{S}(P_{\nu},\mathcal{I}_{\nu},t),\]
where $r_{\nu,\mu}$ is the maximal degree $\lfloor \frac{\dim\nu - \dim\mu -1}{2} \rfloor$ of $g([\mu,\nu],t)$ and where $\tau_{=i}$ selects the term of degree $i$. By induction, this becomes
\[ \tau_{>(d-\dim\mu)/2}\, h_{P_{\mu}}^*(t) - \sum_{\substack{\nu \in \mathcal{P}_f  \\ \nu > \mu }}  \sum_{i=0}^{r_{\nu,\mu}} \tau_{=i}\,g([\mu,\nu],t) \cdot\tau_{> (d-\dim\mu)/2 -i}\, \widetilde{S}(P_{\nu},t)\]
\[ = \tau_{>(d-\dim\mu)/2} \biggl(   h_{P_{\mu}}^*(t) - \sum_{\substack{\nu \in \mathcal{P}_f  \\ \nu > \mu }}  g([\mu,\nu],t) \, \widetilde{S}(P_{\nu},t) \biggr). \tag{6}\] 
Note that $P_{\emptyset}$ has no vertices in the coordinate hyperplanes apart from the origin and the standard basis vectors $\mathbf{e}_i$, since $f$ is convenient. Hence all nonempty faces of $P_{\emptyset}$ that are contained in the order ideal $star_{\partial P_{\emptyset}}(\mathbf{0})$ are unimodular simplices and hence their $\widetilde{S}$-polynomial is zero. So (6) equals
\[  \tau_{>(d-\dim\mu)/2} \biggl(   h_{P_{\mu}}^*(t) - \sum_{\emptyset < F < P_{\mu}}   g([F,P_{\mu}]^*,t) \, \widetilde{S}(F,t) \biggr). \tag{7}\]
Moreover, $ \deg g([\emptyset,P_{\mu}]^*,t) \leq (d-\dim\mu)/2$ and hence (7) equals
\[ \tau_{>(d-\dim\mu)/2} \biggl(   h_{P_{\mu}}^*(t) - \sum_{\emptyset \leq F < P_{\mu}}   g([F,P_{\mu}]^*,t) \, \widetilde{S}(F,t) \biggr)\]
But by Proposition 2.16 again, this is $\tau_{>(d-\dim\mu)/2}\,\widetilde{S}(P_{\mu},t)$ and this has nonnegative coefficients.
\item Secondly, assume that $\Gamma(f)$ has a unique maximal compact face $\nu$. Then for all $\mu\in \mathcal{P}_f\cup \{\emptyset\}$ we have
\[ \widetilde{S}(P_{\mu},star_{\partial P_{\emptyset}}(\mathbf{0})\cap \mathcal{P}(P_{\mu}),t)  =  \widetilde{S}(P_{\mu},P_{\nu},t)    \] and Corollary 2.15 shows immediately that this has nonnegative coefficients.  
\end{enumerate}
This ends the proof of Theorem 5.3. \hfill $\blacksquare$\\

\noindent \textbf{5.5. Example.} The above proof shows that if $\Gamma(f)$ has a unique compact face, then \textsl{all} signs of the local contribution $E_{st,\mathbf{0}}(X_f;u,v)=\sum_{i,j} c_{i,j} u^iv^j$ are `right' in the sense that $(-1)^{i+j}c_{i,j}\geq 0$ for all $i,j$. Now we give a concrete example where the local contribution of
a convenient nondegenerate singularity with crepant first Varchenko subdivision does have `wrong' signs in low degree. Put
\[ f = x_1^2 + x_2^{12} + x_3^{12} +x_4^{12} +x_5^{12} + x_6^{12} +(x_5x_6)^3  \in \mathbb{C}[x_1,\ldots,x_6].\]
The Newton polyhedron $\Gamma(f)$ has two compact facets (i.e.\ faces of codimension 1). One of them has vertices
coming from the monomials
\[ x_1^2,x_2^{12},x_3^{12},x_4^{12},x_5^{12},(x_5x_6)^3       \]
and for the other one just replace $x_5^{12}$ by $x_6^{12}$. In particular, considered as lattice polytopes all
compact faces are simplices. The 1-dimensional cones of the first Varchenko subdivision
$\Delta_f$ are generated by a standard basis vector or by
\[  (6,1,1,1,1,3) \text{ or } (6,1,1,1,3,1).\]
It is easy to check that $\Delta_f$ is indeed crepant. To compute the local contribution
$E_{st,\mathbf{0}}(X_f;u,v)$ one can use the formula from Theorem 4.14 together with the recursion formula from Proposition
2.16. Note that all the involved $g$-polynomials are 1 since the compact faces of the Newton polyhedron are all
simplices (2.5 (2)). To compute the necessary $h^*$-polynomials we used the computer program \textit{Normaliz 2.0} by
Bruns and Ichim \cite{normaliz}. After a rather long computation one finds
\begin{eqnarray*}
E_{st,\mathbf{0}}(X_f;u,v) &=& 3(uv)^4+4(uv)^3-111u^3v^2-111u^2v^3+750u^3v\\ & & +3495(uv)^2+750uv^3
\mathbf{+111u^2v+111uv^2}+3uv+1
\end{eqnarray*}
and thus the sign of $u^iv^j$ is not equal to $(-1)^{i+j}$ for $i+j=3$. The idea of this example is the following, using
notations as above. Consider the compact face $\tau$ equal to the intersection of the two compact facets. Its
associated polytope $P_{\tau}$ has $(6,1,1,1,1,3)$ and $(6,1,1,1,3,1)$ as vertices. It is 1-dimensional and it
has one lattice point in its relative interior. So its $h^*$-polynomial equals $1+t$ and hence
$\widetilde{S}(P_{\tau},\mathcal{I}_{\tau},uv)=uv-1$. This $-1$ gives the wrong signs. It has such a big
influence because there are only 5 nonzero terms in the sum of the formula of Theorem 4.14 and the 4
others come from odd-dimensional faces (including the empty set).

\newpage

\noindent \textbf{5.5. Conclusive remarks.} 
\begin{itemize} 
\item[(1)] I think that Theorem 5.3 is valid for all isolated nondegenerate singularities with crepant first Varchenko subdivision. To generalise the proof of Theorem 5.3 one would need that \[\tau_{>(d-\dim\mu)/2}\,\widetilde{S}(P_{\mu},star_{\partial P_{\emptyset}}(\mathbf{0})\cap \mathcal{P}(P_{\mu}),t)\] has nonnegative coefficients, which I could only prove in the two cases of the theorem. In general, I would guess that the inequality 
\[\tau_{>(d-\dim\mu)/2}\,\widetilde{S}(P_{\mu},star_{\partial P_{\emptyset}}(\mathbf{0})\cap \mathcal{P}(P_{\mu}),t) \geq \tau_{>(d-\dim\mu)/2}\,\widetilde{S}(P_{\mu},t)        \] holds coefficientwise. Note that we proved an equality here for the first case of the theorem. The inequality holds for the second case of the theorem by Proposition 2.14 and by the nonnegativity of the coefficients of relative $g$-polynomials.
\item[(2)] Note that there is a big similarity between the formula of Theorem 4.14 and the formula for the global stringy $E$-function of a nondegenerate Calabi-Yau hypersurface in the toric variety associated to a reflexive polytope \cite[Thm.\ 7.2]{BorisovMavlyutov}.
It is interesting to compare this also with the combinatorial definition of the stringy $E$-function of a Gorenstein lattice polytope
by Batyrev and Nill \cite[Def.\ 4.8]{BatyrevNill}. It is unclear to me what this similarity might mean or suggest. We mention that Batyrev and Nill also formulate an interesting conjecture about their combinatorial stringy $E$-function \cite[Conj.\ 4.10]{BatyrevNill}\\
\end{itemize}

\footnotesize{

\end{document}